\documentclass[12pt,a4paper]{article}
\usepackage{amssymb}
\usepackage{times} 
\usepackage[all]{xy} 
\newcommand{\n}{\newcommand}
\n{\aee}{\begin{enumerate}} \n{\zee}{\end{enumerate}}
\n{\bb}{\bigskip}
\n{\cl}{\centerline}
\n{\del}{\partial}
\n{\D}{\Delta}
\n{\e}{\equiv}
\n{\ep}{\varepsilon}
\n{\f}{\varphi}
\n{\fr}{\frac}
\n{\h}{\widehat}
\n{\ii}{\hskip1em\relax}
\n{\iso}{\stackrel{\sim}{\rightarrow}}
\n{\m}{\bmod}
\n{\mc}{\mathcal}
\n{\mono}{\rightarrowtail}
\n{\Z}{\mathbb Z}
\n{\Zp}{\Z_{(p)}}

\n{\aeq}{\begin{equation}} \n{\zeq}{\end{equation}}

\newtheorem{thm}{Theorem} \n{\at}{\begin{thm}}\n{\zt}{\end{thm}}
\newtheorem{lem}[thm]{Lemma} \n{\al}{\begin{lem}}\n{\zl}{\end{lem}}
\newtheorem{cor}[thm]{Corollary} \n{\ac}{\begin{cor}}\n{\zc}{\end{cor}}
\newtheorem{prop}[thm]{Proposition} \n{\ap}{\begin{prop}}\n{\zp}{\end{prop}}
\newtheorem{defn}[thm]{Definition} \n{\ad}{\begin{defn}}\n{\zd}{\end{defn}}

\begin{document} 

\cl{\Huge Integral Congruences}\bb\bb

To each $i,j$ belonging to some set of integers, attach the integer $a_{ij}$. 
For the congruences 

$$x_j-x_i\e a_{ij}\m(i,j)$$

to be solvable, we must have

$$a_{ij}+a_{jk}\e a_{ik}\m (i,j,k)$$

for all $i,j,k$.

\at This condition is sufficient.\zt

\ii We can generalize the question as follows. Let $n$ be a positive
integer and  $(a(i_0,\dots,i_n))_{i_0,\dots,i_n \in I}$ a family of
integers. For the congruences 

$$\sum_j\ (-1)^j\ x(i_0,\dots,\h{i_j},\dots,i_n)
\e a(i_0,\dots,i_n)\m (i_0,\dots,i_n)$$

to be solvable, we must have 

$$\sum_j\ (-1)^j\ a(i_0,\dots,\h{i_j},\dots,i_{n+1})
\e 0\m (i_0,\dots,i_{n+1})$$

for all $i_0,\dots,i_{n+1}\in I$. 

\at This condition is sufficient.\label{Hn}\zt

\ii To push the generalization further, let $A$ be an abelian group,
$(A_i)_{i\in I}$ a family of subgroups, $n$ a positive integer and 
$(a(i_0,\dots,i_n))_{i_0,\dots,i_n \in I}$ a family of elements of $A$. 
For the congruences 

$$\sum_j\ (-1)^j\ x(i_0,\dots,\h{i_j},\dots,i_n)
\e a(i_0,\dots,i_n)\m A_{i_0}+\cdots+A_{i_n}$$

to be solvable, we must have 

$$\sum_j\ (-1)^j\ a(i_0,\dots,\h{i_j},\dots,i_{n+1})
\e 0\m A_{i_0}+\cdots+A_{i_{n+1}}$$

for all $i_0,\dots,i_{n+1}\in I$. 

\at\label{3} This condition is {\em not} sufficient in general. \zt

\ii The above statement and its proof were suggested by Pierre Schapira and
Terence Tao. --- Recall that a lattice is {\bf distributive} if one of the two
operations is distributive over the other, or, equivalently, if each
operation is distributive over the other. 

\at\label{4} If the indexing set $I$ is finite and if the $A_i$ generate a
 distributive lattice, then the above condition is sufficient. 
\zt

\ii The above statement and its proof were suggested by Anton Deimtar. 

\section{Refinements}\label{r}

Let $A$ be an abelian group and $(A_i)_{i\in I}$ a family of subgroups. 
Form the cochain complex 

$$C^n(I):=\prod_{i_0,\dots,i_n\in I}\ \fr{A}{A_{i_0}+\cdots+A_{i_n}}\quad,\quad
d:C^{n-1}(I)\to C^n(I),$$

$$(df)(i_0,\dots,i_n)\e\sum_j\ (-1)^j\ f(i_0,\dots,\h{i_j},\dots,i_n)
\m A_{i_0}+\cdots+A_{i_n}.$$

Let $(B_j)_{j\in J}$ be another family of subgroups. A map $\tau$ from $J$
to $I$ is a {\bf refinement map} if 

$$A_{\tau j}\subset B_j\quad\forall\ j\in J,$$

and $J$ is a {\bf refinement} of $I$ if there is such a refinement map. A
refinement map $\tau$ induces a cochain map, still denoted $\tau$, from
$C^n(I)$ to $C^n(J)$ defined by   

$$(\tau f)(j_0,\dots,j_n)\e f(\tau j_0,\dots,\tau j_n)
\m B_{j_0}+\cdots+B_{j_n},$$

and thus a morphism $\tau^*$ from $H^n(I)$ to $H^n(J)$. \bb

\ii Let $\sigma$ be another refinement map from $J$ to $I$. We claim
$\sigma^*=\tau^*$. Define the morphism $h$ from $C^n(I)$ to
$C^{n-1}(J)$ by  

$$(hf)(j_0,\dots,j_{n-1}):=$$

$$\sum_k\ (-1)^k\ f(\tau j_0,\dots,\tau j_k,\sigma j_k,\dots,\sigma j_{n-1})
\m  B_{j_0}+\cdots+B_{j_{n-1}}.$$ 

We have

$$dh+hd=\sigma-\tau,$$

which implies the claim. \bb

\ii Let $(A_i)_{i\in I}$ be as above, put

$$J:=\{A_i\ |\ i\in I\},$$

let $(B_j)_{j\in J}$ be the tautological family, let $\pi$ be the natural
surjection from $I$ onto $J$ and $\iota$ a section of $\pi$. Then $\pi$ and
$\iota$ are refining maps, and $C(I)$ and $C(J)$ are homotopy equivalent. 

\section{Coefficient Systems}

Recall the a {\bf simplicial complex} is a set $K$ equipped with a set of
nonempty finite subsets, called {\bf simplices} subject to the condition
that a nonempty subset of a simplex is a simplex. A map between two
simplicial complexes is {\bf simplicial} if it maps simplices to
simplices. For $n\ge0$ let $\D_n$ be the set $\{0,1,\dots,n\}$ equipped
with the simplicial structure giving the status of simplex to all nonempty
finite subsets, say that a a {\bf singular $n$-simplex} of $K$ is a
simplicial map from $\D_n$ to $K$, and let $S_n(K)$ be the set of  singular
$n$-simplices of $K$. For $n>0$ and $i\in\D_n$ let $f_i$
be the increasing map from $\D_{n-1}$ to $\D_n$ missing the vertex $i$. \bb

\ii Let $C_n(K)$ be the free abelian group generated $S_n(K)$ and consider
the morphisms

$$\del_n=\del:C_n(K)\to C_{n-1}(K),\quad 
\del s:=\sum_i\ (-1)^i\ sf_i\ \forall\ s\in K_n$$

for $n>0$, and 

$$\ep:C_0\to \Z,\quad s\mapsto1\ \forall\ s\in S_0(K).$$

Then $(C_*(K),\del,\ep)$ is an augmented chain complex. 

\al Let $s$ be in $S_n(K)$ and $\Phi(s)$ the subgroup of $C_*(K)$ generated
 by the singular simplices of the form $sf$ with $f:\D_k\to\D_n$. Then
 $\Phi(s)$ is an acyclic augmented subcomplex of $C_*(K)$.
\zl

{\bf Proof.} It is clear that $\Phi(s)$ is a subcomplex of $C_*(K)$, and
that

$$sf\mapsto(s0,sf0,\dots,sfn)$$

is a homotopy from the identity of $\Phi(s)$ to 0. QED\bb

\ii Let $\le$ be an ordering on $K$ and $\f$ be the $\Z$-linear endomorphism
of $C_*(K)$ defined by\bb

\ii $\f s:=s$ if $s$ is noninjective,\bb

\ii $\f(s0,\dots,sn):=s-(-1)^\sigma\ (s\sigma0,\dots,s\sigma n)$ if
$s$ is injective, $\sigma$ is the permutation characterized by
$s\sigma0<\cdots<s\sigma n$, and $(-1)^\sigma$ is the signature of
$\sigma$. \bb

\ii Then $\f$ is an endomorphism of the augmented complex
$(C_*(K),\del,\ep)$. 

\al There is a homotopy $h$ from $\f$ to 0 such that $hs$ is in $\Phi(s)$ 
 for all $s\in S_n(K)$.
\zl

{\bf Proof.} Put $C_n:=C_n(K)$, let $h_0$ be the zero morphism from $C_0$ to
$C_1$, and assume that we have morphisms $h_n:C_n\to C_{n+1}$ satisfying 

$$\del_{n+1}\ h_n=\f_n-h_{n-1}\ \del_n$$

for $n<k$. We want to define $h_k:C_k\to C_{k+1}$ in such a way that we
have 

$$\del_{k+1}\ h_k=\f_k-h_{k-1}\ \del_k.$$

Observe

\begin{eqnarray*}
 \del_k\ (\f_k-h_{k-1}\ \del_k)&=&\f_{k-1}\ \del_k-\del_k\ h_{k-1}\
  \del_k\\
 &=&\f_{k-1}\ \del_k-(\f_{k-1}-h_{k-2}\ \del_{k-1})\ \del_k\\
 &=&0
\end{eqnarray*}

and use the previous Lemma. QED\bb

\ii Let $A$ be  a ring (commutative with 1). Define the category $K^\D$ as
follows. The objects of $K^\D$ are the singular simplices of $K$. The
morphisms from $s\in S_n(K)$ to $t\in S_k(K)$ are the maps $f$ from $\D_k$
to $\D_n$ such that $t=sf$, the composition being the obvious one. Say that a 
{\bf coefficient system} over $K$ is a functor from $K^\D$ to the category
of $A$-modules. \bb 

\ii If $V$ is a coefficient system over $K$, then we denote by $V(s)$ the
$A$-module attached to $s\in S_n(K)$ and by $V(s,f)$ the morphism from
$V(s)$ to $V(sf)$ associated with the map $f$ from $\D_k$ to $\D_n$. \bb

\ii Denote by $c=(c(s))_{s\in S_n(K)}\in C^n(K,V)$ the vectors of the
$A$-module  

$$C^n(K,V):=\prod_{s\in S_n(K)}V(s).$$

and consider the morphisms

$$d_{n-1}=d:C^{n-1}(K,V)\to C^n(K,V),\quad
(dc)(s):=\sum_i\ (-1)^i\ V(s,f_i)\ c(s)$$

for $n>0$. Then $(C^*(K,V),d)$ is a cochain complex. \bb

\ii Assume $V_{s\sigma}=V_s$ whenever $s$ is injective and $\sigma$ is a
permutation. Say that the cochain $c$ of $C^n(K,V)$ is {\bf alternating}
if\bb 

\ii (a) $c(s)=0$ whenever $s$ is noninjective,\bb

\ii (b) $c(s\sigma)=-c(s)$ whenever $s$ is injective and $\sigma$ is an odd
permutation. \bb

The alternating cochains form a subcomplex $C'^*(K,V)$ of $C^*(K,V)$. \bb

\ii If $\psi$ is a $\Z$-linear map from $C_p(K)$ to $C_q(K)$ of the form

$$\psi s=\sum_{f:\D_q\to\D_p}\lambda_{s,f}\ sf$$

with $\lambda_{s,f}\in\Z$, if $c$ is in $C^q(K,V)$ and $s$ in $S_p(K)$, we
put 

$$(h'c)(s):=\sum_f\ \lambda_{s,f}\ V(s,f)\ c(s).$$

\ii We clearly have

\ap In the above notation $1-\f'$ is a projector onto $C'^*(K,V)$. Moreover
if $h$ is the homotopy of the previous Lemma, then $h'$ is a homotopy from
$1-\f'$ to the identity. In particular the inclusion of $C'^*(K,V)$ into
$C^*(K,V)$ is a quasi-isomorphism.\label{alt}
\zp

\section{Rings}

By ``ring'' we mean``commutative ring with 1''. A {\bf domain} is a nonzero ring
which has no nontrivial zero divisors. An ideal of a ring is {\bf prime} if
the quotient is a domain. 

\al The inverse image of a prime ideal under a ring morphism is
 prime.\label{ii}
\zl 

{\bf Proof.} Left to the reader.\bb

\ii A subset $S$ of a ring $A$ is {\bf multiplicative} if it contains 1 and is
closed under multiplication.  Let $A$ be a ring, $S$ a multiplicative
subset and $V$ an $A$-module. For $(v_1,s_1),(v_2,s_2)\in V\times S$ write
$(v_1,s_1)\sim(v_2,s_2)$ if there is an $s\in S$ such that
$s(s_1v_2-s_2v_1)=0$. This is an equivalence relation. Let $S^{-1}V$ be the
quotient and denote the class of $(v,s)$ by $v/s$, by $\fr{v}{s}$ or by
$s^{-1}v$. For $(v_1,s_1),(v_2,s_2)\in V\times S,a\in A,v\in V$ put

$$\fr{v_1}{s_1}+\fr{v_2}{s_2}=\fr{s_2v_1+s_1v_2}{s_1s_2}\quad,\quad
\fr{a}{s_1}\ \fr{v}{s_2}=\fr{av}{s_1s_2}\quad.$$

These formulas equip respectively $S^{-1}A$ and $S^{-1}V$ with well defined
structures of ring and $(S^{-1}A)$-module. We say that $S^{-1}A$ and
$S^{-1}A$ are respectively the 
{\bf ring of fractions of $A$ with denominators in} $S$ and the 
{\bf module of fractions of $V$ with denominators in} $S$. The map 

$$i_V:V\to S^{-1}V,\quad v\mapsto\fr{v}{1}$$

is called the {\bf canonical morphism}.

\al The canonical morphism $i_A$ is a ring morphism and $i_V$ \label{ker}
 is an $A$-module morphism whose kernel consists of those elements of $V$
 which have an annihilator in $S$. 
\zl 

{\bf Proof.} Left to the reader.

\al[Universal Property]\label{up}
 Any ring morphism $A\to B$ mapping $S$ into the multiplicative group
 $B^\times$ factors uniquely through $S^{-1}A$. Under the same assumption, if
 $V$ is an $A$-module and $W$ a $B$-module, then any $A$-module morphism  
 $V\to W$ factors uniquely through $S^{-1}V$, the induced map $S^{-1}V\to W$
 being an $(S^{-1}A)$-module morphism. 
\zl

{\bf Proof.} Left to the reader.

\al The $(S^{-1}A)$-modules $S^{-1}V$ and $S^{-1}A\otimes_AV$ are
      canonically isomorphic. 
\zl

{\bf Proof.} Left to the reader.\bb

\ii Recall that an $A$-module $F$ is {\bf flat} if for all $A$-module
monomorphism $W\mono V$ the obvious morphism 
$F\otimes_A W\to F\otimes_A V$ is injective. The Lemma below, which follows
from the previous one, will be freely used. 

\al The $A$-module $S^{-1}A$ is flat.\label{f}\zl

\ii Let $A$ be a ring, let $S$ be a multiplicative subset, let $i$ be the canonical
morphism $i_A$, let $\mc{J}$ be the set of all ideals of $S^{-1}A$
and $\mc{I}$ the set of those ideals $I$ of $A$ such that the nonzero
elements of $A/I$ have no annihilators in $S$. 

\al We have\label{90} 
 \aee
  \item $S^{-1}I=(S^{-1}A)i(I)$ for all ideal $I$ of $A$;
  \item $S^{-1}i^{-1}(J)=J$ for all $J\in\mc{J}$;
  \item if $I$ is an ideal of $A$, the ideal $i^{-1}(S^{-1}I)$ of $A$ is
    formed by those $a$ in $A$ whose image in $A/I$ has an annihilator in
    $S$; 
  \item $i^{-1}$ maps $ \mc{J}$ bijectively onto $\mc{I}$ and the
    inverse bijection is given by $S^{-1}$:
    $$\xymatrix{\mc{I}\ar@<1ex>[r]^{S^{-1}}&\mc{J}.\ar@<1ex>[l]^{i^{-1}}}$$
  \item for $J\in\mc{J}$ and $I:=i^{-1}(J)$, the canonical morphism
    from $S^{-1}A$ to $S^{-1}(A/I)$ is surjective, its kernel is $J$ and
    $i_{A/I}$ is injective; 
  \item a prime ideal of $A$ is in $\mc{I}$ iff it is disjoint from $S$;
  \item $S^{-1}$ and $i^{-1}$ induce inverse bijections between the set of
    prime ideals of $A$ disjoint from $S$ and the set of prime ideals of
    $S^{-1}A$. \label{7}
 \zee
\zl

{\bf Proof.} Follows from the previous Lemmas. QED

\al An ideal of $A$ is prime iff its complement is multiplicative.\label{pr}\zl 

{\bf Proof.} Left to the reader.\bb

\ii If $S\subset A$ is the complement of a prime ideal $P$ and $V$ an
$A$-module, we put 

\aeq\label{AP}A_P:=S^{-1}A,\quad V_P:=S^{-1}V.\zeq

\al\label{p} Let $p$ be a prime number and $(F_i)$ a family of finite
abelian groups. Then the natural morphism
$$\left(\prod_iF_i\right)_{(p)}\to\prod_i\ (F_i)_{(p)}$$
is bijective.
\zl

{\bf Proof.} We can assume either that each $F_i$ is a $p$-group, or
that each $F_i$ is a $p_i$-group for some prime $p_i\not=p$. QED 

\section{The local case}

Let $p$ be a prime number. Part~\ref{7} of Lemma~\ref{90} implies that the
ideals of $\Zp$ are precisely the powers of $p\Zp$. Let $J_1,J_2,\dots$ be
a (possibly finite) decreasing sequence of ideals of $\Zp$, and form the
cochain complex  

$$C^n:=\prod_{i_0<\cdots<i_n\in I}\ \fr{\Zp}{J_{i_0}}\quad,\quad
d:C^{n-1}\to C^n,$$

$$(df)(i_0,\dots,i_n)\e\sum_j\ (-1)^j\ f(i_0,\dots,\h{i_j},\dots,i_n)\m
J_{i_0}.$$

\ap We have $H^n(C)=0$ for $n>0$.\zp

{\bf Proof.} We denote by the same symbol an element of $\Zp$ and its class
modulo $J_i$.  Let $a$ be an $n$-cocycle of $C$, and define the element
$x(i_1,\dots,i_n)$ of $\Zp$ inductively on $i_0$ as follows. Put first

$$x(1,i_2,\dots,i_n):=0.$$

For $i_1>1$ set $i_0:=i_1-1$ and

$$x(i_1,\dots,i_n):=
-\sum_{j=1}^n\ (-1)^j\ x(i_0,\dots,\h{i_j},\dots,i_n)+a(i_0,\dots,i_n).$$

Given $1<i_1<\cdots<i_n$ prove

$$x(i_1,\dots,i_n)\e
-\sum_{j=1}^n\ (-1)^j\ x(i_0,\dots,\h{i_j},\dots,i_n)+a(i_0,\dots,i_n)
\m J_{i_0}$$

successively for $i_0=i_1-1,i_1-2,\dots,1$. QED

\section{Proof of the Theorems}

{\bf Proof of Theorem \ref{Hn}.}
Let $C$ be the cochain complex obtained by replacing $A$ by $\Z$ and $A_i$ 
by $(i)$ in the definition of the cochain complex $C(I)$ at the outset of
Section~\ref{r}. Let $n$ be a positive integer. By Lemma \ref{f} we have
$H^n(C)_{(p)}=H^n(C_{(p)})$. By Lemma~\ref{ker} an abelian group $A$ is
trivial iff $A_{(p)}=0$ for all $p$. The triviality of $H^n(C_{(p)})$
follows from Section~\ref{r}, Proposition~\ref{alt}, Part~\ref{7} of
Lemma~\ref{90}, Lemma~\ref{p}, and the above Proposition. \bb

{\bf Proof of Theorem \ref{3}.}
Let $V$ be a three dimensional vector space; let $L_1$, $L_2$, $L_3$, $L_4$
be four lines in general position; let $C$ be the cochain complex denoted
$C(I)$ in the beginning of Section~\ref{r}. We claim $H^1(C)\not=0$. By
Proposition~\ref{alt}  we can replace the condition $i_0,\dots,i_n\in I$ in
the definition of $C$ by the condition $i_0<\cdots<i_n\in I$; in
particular we get $C^n=0$ for $n>1$, and we claim that the coboundary 
$d$ from  

$$C^0=\fr{V}{L_1}\oplus\fr{V}{L_2}\oplus\fr{V}{L_3}\oplus\fr{V}{L_4}$$

to 

$$C^1=\fr{V}{L_1+L_2}\oplus\fr{V}{L_1+L_3}\oplus\fr{V}{L_1+L_4}
\oplus\fr{V}{L_2+L_3}\oplus \fr{V}{L_2+L_4}\oplus\fr{V}{L_3+L_4}$$\bb

is {\bf not} onto. But his follows from the fact that the kernel of $d$ contains
$V$. \bb

{\bf Proof of Theorem \ref{4}.}
As a general notation, write $H^n((A_i),A)$ for the $n$-coho\-molo\-gy of
the cochain complex defined at the beginning of Section~\ref{r}. Theorem
\ref{4} results from Theorem~\ref{Hn} and the following fact, whose proof is
left to the reader. 

\al In the notation of the beginning of Section~\ref{r}, assume that the
lattice generated by the $A_i$ is distributive, let $B$ be one of the $A_i$,
and $n$ a positive integer. Then 
$$H^n((B\cap A_i),A\big)=0
=H^n\Bigg(\bigg(\fr{B+A_i}{B}\bigg),\fr{A}{B}\Bigg)
\Rightarrow H^n((A_i),A)=0.$$
\zl


\vfill\hfill Pierre-Yves Gaillard, Universit\'e Nancy 1, France

\end{document}